\documentclass{ifacconf}

\usepackage{graphicx}      % include this line if your document contains figures
\usepackage{natbib}        % required for bibliography
\usepackage{amsmath}
\usepackage{amsfonts}
\begin{document}
\begin{frontmatter}

\title{Optimal Control of Connected and Automated Vehicles at Roundabouts: An Investigation in a Mixed-Traffic Environment\thanksref{footnoteinfo}}

\thanks[footnoteinfo]{This research was supported by US Department of Energy's
SMART Mobility Initiative.}

\author[First]{Liuhui Zhao} 
\author[First]{Andreas Malikopoulos} 
\author[Second]{Jackeline Rios-Torres}

\address[First]{University of Delaware, Newark, DE 19716 USA (e-mails: lhzhao@udel.edu; andreas@udel.edu).}
\address[Second]{Oak Ridge National Laboratory, Oak Ridge, TN 37830 USA (e-mail: riostorresj@ornl.gov)}

\begin{abstract}                % Abstract of not more than 250 words.
Connectivity and automation in vehicles provide the most intriguing opportunity for enabling users to better monitor transportation network conditions and make better operating decisions to improve safety and reduce pollution, energy consumption, and travel delays. This paper investigates the implications of optimally coordinating vehicles that are wirelessly connected to each other in roundabouts to achieve a smooth traffic flow. We apply an optimization framework and an analytical solution that allows optimal coordination of vehicles for merging in such traffic scenario. The effectiveness of the proposed approach is validated through simulation and it is shown that fully coordinated vehicles reduce total travel time by 51\% and fuel consumption by 35\%. Furthermore, we study the influence of vehicle coordination in a mixed-traffic environment and compare the network performance under different market penetration rates of connected and automated vehicles (CAVs). For this particular study with near-capacity demand, due to extremely unstable traffic, the results show that even with high penetration of CAVs (e.g., 80\%), travel time and fuel savings are much less than a network of CAVs.
\end{abstract}


\end{frontmatter}
%===============================================================================

\section{Introduction}
We are currently witnessing an increasing integration of energy and transportation, which, coupled with human interactions, is giving rise to a new level of complexity in the next generation transportation systems; see \cite{Malikopoulos2015}. The common thread that characterizes energy efficient mobility systems is their interconnectivity which enables the exchange of massive amounts of data; this, in turn, provides the opportunity for a novel computational framework to process such information and deliver real-time control actions that optimize energy consumption and associated benefits. As we move to increasingly complex transportation systems new control approaches are needed to optimize the impact on system behavior of the interplay between vehicles at different traffic scenarios; see \cite{Malikopoulos2017}. 

Intersections, roundabouts, merging roadways, speed reduction zones are the primary sources of bottlenecks that contribute to traffic congestion created by the drivers’ responses; see \cite{Rios-Torres2017, Malikopoulos2013}. In 2015, congestion caused people in urban areas in US to spend 6.9 billion hours more on the road and to purchase an extra 3.1 billion gallons of fuel, resulting in a total cost estimated at \$160 billion; see \cite{Schrank2015}. 

Roundabouts generally provide better operational and safety characteristics over other typical intersections; see \cite{Flannery1997, Flannery1998, Al-Madani2003, Sisiopiku2001}. However, the increase of traffic may become a concern for roundabouts due to their geometry and priority system - even with moderate demands, some roundabouts can still reach capacity; see \cite{Hummer2014, Yang2004}. Moreover, traffic from minor-road approaches may experience significant delay if the circulating flow is heavy. Previous research has focused mainly on enhancing roundabout mobility and safety with improved metering, or traffic signal controls, e.g., \cite{Hummer2014, Yang2004, Martin-Gasulla2016, Xu2016}. These efforts investigate the potential of metering signals in improving roundabout operations during rush hours. \cite{Hummer2014} tested a metering approach for a single-lane and a two-lane roundabout models with different levels of approaching traffic demand. \cite{Martin-Gasulla2016} studied the benefits of metering signals for roundabouts with unbalanced flow patterns. \cite{Yang2004} proposed a traffic-signal control algorithm to eliminate the conflict points and weaving sections for multi-lane roundabouts by introducing a second stop line for left-turn traffic. \cite{Xu2016} suggested a multi-level control system that combines metering signalization with full actuated control to serve different time periods throughout the day.

Connected and automated vehicles (CAVs) provide the most intriguing and promising opportunity for enabling users to better monitor transportation network conditions and make better operating decisions to reduce energy consumption, greenhouse gas emissions, travel delays and improve safety. Given the recent technological developments, several research efforts have considered approaches to achieve safe and efficient coordination of merging maneuvers with the intention to avoid severe stop-and-go driving. One of the very early efforts in this direction was proposed in 1969 by \cite{Athans1969}. Assuming a given merging sequence, Athans formulated the merging problem as a linear optimal regulator, proposed by \cite{Levine1966}, to control a single string of vehicles, with the aim of minimizing the speed errors that will affect the desired headway between each consecutive pair of vehicles. Later, \cite{Schmidt1983} proposed a two-layer control scheme based on heuristic rules that were derived from observations of the non-linear system dynamics behavior. Similar to the approach proposed by \cite{Athans1969}, \cite{Awal2013} developed an algorithm that starts by computing the optimal merging sequence to achieve reduced merging times for a group of vehicles that are closer to the merging point.

Several research efforts have been reported in the literature proposing different approaches on coordinating CAVs at different transportation segments, e.g., intersections, roundabouts, merging roadways, speed reduction zones, with the intention to improve traffic flow. \cite{Kachroo1997} proposed a longitudinal and lateral controller to guide the vehicle until the merging maneuver is completed. Other efforts have focused on developing a hybrid control aimed at keeping a safe headway between the vehicles in the merging process, see \cite{Kachroo1997, Antoniotti1997}; or developing three levels of assistance for the merging process to select a safe space for the vehicle to merge; see \cite{Ran1999}. Other authors have explored virtual vehicle platooning for autonomous merging control, e.g., \cite{Lu2000, Dresner2004}, where a controller identifies and interchanges appropriate information between the vehicles involved in the merging maneuver while each vehicle assumes its own control actions to satisfy the assigned time and reference speed. 

There has been a significant amount of work on vehicle coordination at urban intersections. \cite{Dresner2004} proposed the use of the reservation scheme to control a single intersection of two roads with vehicles traveling with similar speed on a single direction on each road, i.e., no turns are allowed. In their approach, each vehicle is treated as a driver agent who requests the reservation of the space-time cells to cross the intersection at a particular time interval defined from the estimated arrival time to the intersection. Once the reservation system receives the request, it accepts if there is no conflict with the already accepted reservations; otherwise, the request is to be rejected. In case of rejection, the driver agent is required to decelerate and send a new reservation request. Since then, numerous approaches have been reported in the literature, e.g., \cite{DeLaFortelle2010, Dresner2008, Huang2012},  to achieve safe and efficient control of traffic through intersections including extensions of the reservation scheme in \cite{Dresner2004}. Other research efforts have focused on coordinating vehicles at intersections to improve travel time, e.g., \cite{Zohdy2012, Zhu2015, Yan2009}.

In this paper, we address the problem of optimally coordinating a number of CAVs entering a roundabout, so as to improve traffic flow. The objectives of this study are to: 1) formulate the problem of optimal CAV coordination at roundabouts and 2) investigate the implications on fuel consumption and travel time at different market penetration levels.

The contributions of this paper are 1) the formulation of the problem of controlling CAVs before they enter a roundabout, 2) implementing an analytical solution that yields the optimal acceleration/deceleration for each vehicle, and 3) investigating the impact of the optimal solution through simulation under different traffic conditions.

\section{Problem Formulation}

\begin{figure} [!ht]
\begin{center}
\includegraphics[width=8.4cm]{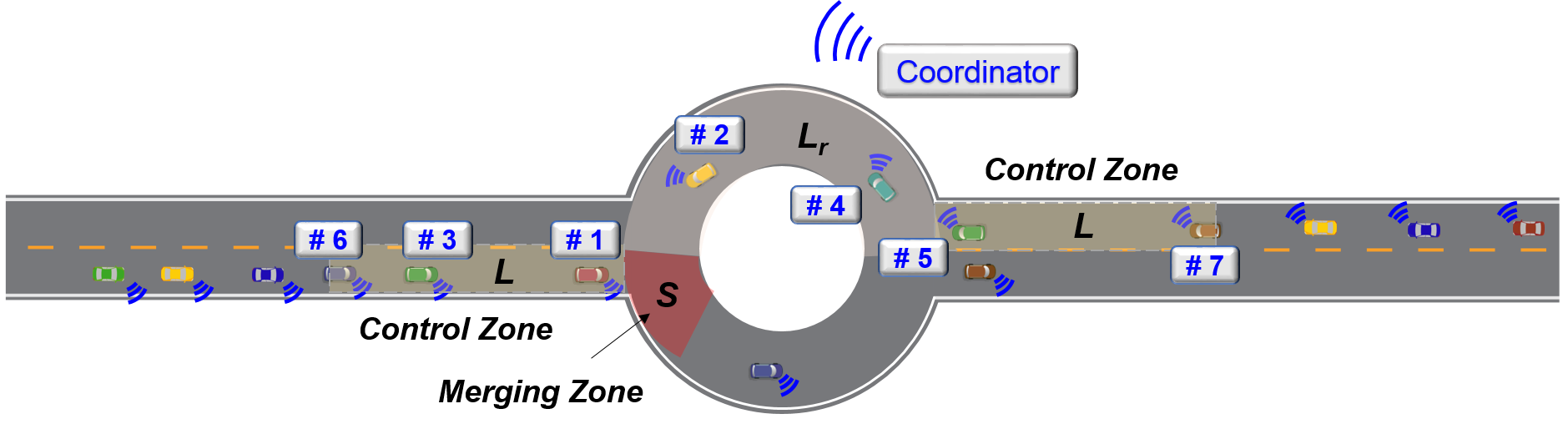}    % The printed column width is 8.4 cm.
\caption{A roundabout with a traffic flow of CAVs.} 
\label{fig:roundaboutcon}
\end{center}
\end{figure}

We consider a single-lane roundabout (Fig.~\ref{fig:roundaboutcon}) where traffic from two different freeways entering the roundabout with a higher speed than the imposed roundabout speed limit. Before the entry of the roundabout, there is a \emph{control zone} and a coordinator that can communicate with the vehicles traveling inside the control zone. Note that the coordinator is not involved in any decision on the vehicle. The region at the roundabout where a potential lateral collision of the vehicles can occur is called \emph{merging zone}. The arc length of the merging zone is $S$, and the length of the control zone is $L$. The arc length from the exit of the control zone to the entry of the merging zone is $L_r$. Note that $L$ could be in the order of hundreds of $m$ depending on the coordinator's communication range capability. 

We consider a number of CAVs $N(t)\in \mathbb{N}$, where $ t\in \mathbb{R}^+$ is the time, entering the control zone. Let $\mathcal{N}(t)={1,\dots,N(t)}$ be a first-in-first-out (FIFO) queue in the control zone. The dynamics of each vehicle $i, i\in \mathcal{N}(t)$, are represented with a state equation

\begin{equation} \label{eq:state}
\dot{x}(t) = f(t, x_i, u_i), ~ x_i(t_i^0) = x_i^0,
\end{equation}

where $x_i(t), u_i(t)$ are the state of the vehicle and control input, $t_i^0$ is the time that vehicle $i$ enters the control zone, and $x_i^0$ is the value of the initial state. For simplicity, we model each vehicle as a double integrator, i.e., $\dot{p}_i = v_i(t)$ and $\dot{v}_i = u_i(t)$, where $p_i(t) \in \mathcal{P}_i, v_i(t) \in \mathcal{V}_i$, and $u_i(t) \in \mathcal{U}_i$ denote the position, speed, and acceleration/deceleration (control input) of each vehicle $i$. Let $x_i(t)=[p_i(t)~ v_i(t)]^T$ denotes the state of each vehicle $i$, with initial value $x_i^0(t)=[0 ~ v_i^0(t)]^T$, taking values in the state space $\mathcal{X}_i=\mathcal{P}_i\times  \mathcal{V}_i$. The sets $\mathcal{P}_i, \mathcal{V}_i$, and $\mathcal{U}_i, i\in \mathcal{N}(t)$, are complete and totally bounded subsets of $\mathbb{R}$. The state space $\mathcal{X}_i$ for each vehicle $i$ is closed with respect to the induced topology on $\mathcal{P}_i\times  \mathcal{V}_i$ and thus, it is compact.

To ensure that the control input and vehicle speed are within a give admissible range, the following constraints are imposed.
\begin{equation}\label{eq:constraints}
\begin{aligned} 
u_{min} &\leq u_i(t) \leq u_{max},~ \text{and} \\
0 &\leq v_{min} \leq v_i(t) \leq v_{max}, ~ \forall t \in [t_i^0, ~ t_i^z]
\end{aligned} 
\end{equation}

where $u_{min}, u_{max}$ are the minimum deceleration and maximum acceleration respectively, and $v_{min}, v_{max}$ are the minimum and maximum speed limits respectively. $t_i^0$ is the time that vehicle $i$ enters the control zone, $t_i^z$ is the time that vehicle $i$ exits the control zone.

To ensure the absence of rear-end collision of two consecutive vehicles traveling on the same lane, the position of the preceding vehicle should be greater than, or equal to the position of the following vehicle plus a predefined safe distance $\delta_i(t)$, where $\delta_i(t)$ is proportional to the speed of vehicle $i$, $v_i(t)$. Thus, we impose the rear-end safety constraint
\begin{equation} \label{eq:safety}
s_i(t)=p_k(t)-p_i(t)\geq \delta_i(t), \forall t \in  [t_i^0, ~ t_i^f]
\end{equation}

where vehicle $k$ is immediately ahead of $i$ on the same road, and $t_i^f$ is the time that vehicle $i$ exits the merging zone. In the aforementioned modeling framework, we assume that each vehicle cruises inside the roundabout at the imposed speed limit, i.e., $v_i(t)=v_r, \forall t \in [t_i^z, ~ t_i^f]$.

We consider the problem of minimizing the control input (acceleration/deceleration) for each vehicle $i$  from the time $t_i^0$ that the vehicle enters the control zone until the time $t_i^z$ that it exits the control zone under the hard safety constraint to avoid rear-end collision. Thus, we formulate the following optimization problem for each vehicle in the queue $\mathcal{N}(t)$
\begin{equation} \label{eq:min}
\begin{aligned}
&\min_{u_i}\frac{1}{2}\int_{t_i^0}^{t_i^z} u_i^2(t)dt \\
&Subject~to: (\ref{eq:state}), (\ref{eq:constraints}) ~\forall i \in \mathcal{N}(t)
\end{aligned}
\end{equation}

The analytical solution of (\ref{eq:min}) without considering state and control constraints was presented in \cite{Rios-Torres2017a} and \cite{Ntousakis2016} for coordinating online CAVs at highway on-ramps and in \cite{Zhang2016} at two adjacent intersections. The solution of the constrained problem at a single intersection was presented in \cite{Malikopoulos2017a}. To implement the analytical solution of  (\ref{eq:min}), each vehicle $i$ needs to compute the time $t_i^z$ at which it will be exiting the control zone. Thus, we introduce the notion of the coordinator to handle the information between the vehicles as follows. When a vehicle reaches the control zone at some instant $t$, the coordinator assigns a unique identity, which is an integer $i$ representing the order of the vehicle in the FIFO queue inside the control zone. Once vehicle $i$ enters the control zone, it shares its time $t_i^z$ that it will be exiting the control zone. Then the vehicle $i+1$ in the queue computes the time $t_{i+1}^z$ that it will exit the control zone with respect to $t_i^z$. Thus, the time $t_i^m$ that each vehicle will be entering the merging zone can be computed directly from $t_i^z$.

In the situation that westbound traffic will enter the roundabout and form a circulating flow with which eastbound vehicles interact in the merging zone, we define that each vehicle $i\in \mathcal{N}(t)$ belongs to either of two different subsets: 1) $\mathcal{L}(t)$ contains all vehicles traveling westbound, and 2) $\mathcal{C}(t)$ contains all vehicles traveling eastbound. Thus, we set
\begin{equation} \label{eq:tm}
t_i^m=t_i^z+\frac{\lambda_i \cdot L_r}{v_r}
\end{equation}

where
\begin{equation} \label{eq:lambda}
\lambda_i = \begin{cases} 0, & \quad \forall i \in \mathcal{C}(t) \\ 1, & \quad  \forall i \in \mathcal{L}(t)
\end{cases}
\end{equation}

$v_r$ is the imposed speed limit inside the roundabout, and $\lambda_i$ is an indicator corresponding to the traveling approach of vehicle $i$. 

The time $t_i^m$ that vehicle $i$ will be entering the merging zone is restricted by the imposed rear-end collision constraint. To ensure that (\ref{eq:safety}) is satisfied at time $t_i^m$ and that $t_i^m$ can be achieved within the imposed control and speed limits, we impose the following conditions:

1) If vehicles $i-1$ and $i$ are traveling on the same road, then vehicles $i-1$ and $i$ should have the predefined safe distance allowable, denoted by $\delta_i(t)$, by the time vehicle $i$ enters the merging zone, i.e.,

\begin{equation} \label{eq:tm_re_1}
\begin{split}
t_i^m=\max \Big\{ \min \{ 
t_{i-1}^m + \frac{\delta_i(t_i^m)}{v_r}, 
t_i^0+\frac{L}{v_{min}} + \frac{\lambda_i \cdot L_r}{v_r} \}, \\
t_i^0+\frac{L}{\bar{v}_i} + \frac{\lambda_i \cdot L_r}{v_r},
t_i^0+\frac{L}{v_{max}}+\frac{\lambda_i \cdot L_r}{v_r} \Big\}
\end{split}
\end{equation}

2) If vehicle $i-1$ and $i$ are traveling on the different roads, we constrain the merging zone to contain only one vehicle so as to avoid a lateral collision. Therefore, vehicle $i$ is allowed to enter the merging zone only when vehicle $i-1$ exits the merging zone, i.e., 

\begin{equation} \label{eq:tm_re_2}
\begin{split}
t_i^m=\max \Big\{ \min \{ 
t_{i-1}^m + \frac{S}{v_r}, 
t_i^0+\frac{L}{v_{min}} + \frac{\lambda_i \cdot L_r}{v_r} \}, \\
t_i^0+\frac{L}{\bar{v}_i} + \frac{\lambda_i \cdot L_r}{v_r},
t_i^0+\frac{L}{v_{max}}+\frac{\lambda_i \cdot L_r}{v_r} \Big\}
\end{split}
\end{equation}

where $\bar{v}_i$ is the average speed for vehicle $i$ traveling from the entry to the exit of control zone. The recursion is initialized whenever a vehicle enters a control zone, i.e., it is assigned $i=1$. In this case, $t_i^m$ can be externally assigned as the desired exit time of this vehicle whose behavior is unconstrained except for (\ref{eq:state}).

\section{Simulation Framework}
\begin{figure} [!ht]
\begin{center}
\includegraphics[width=8.4cm]{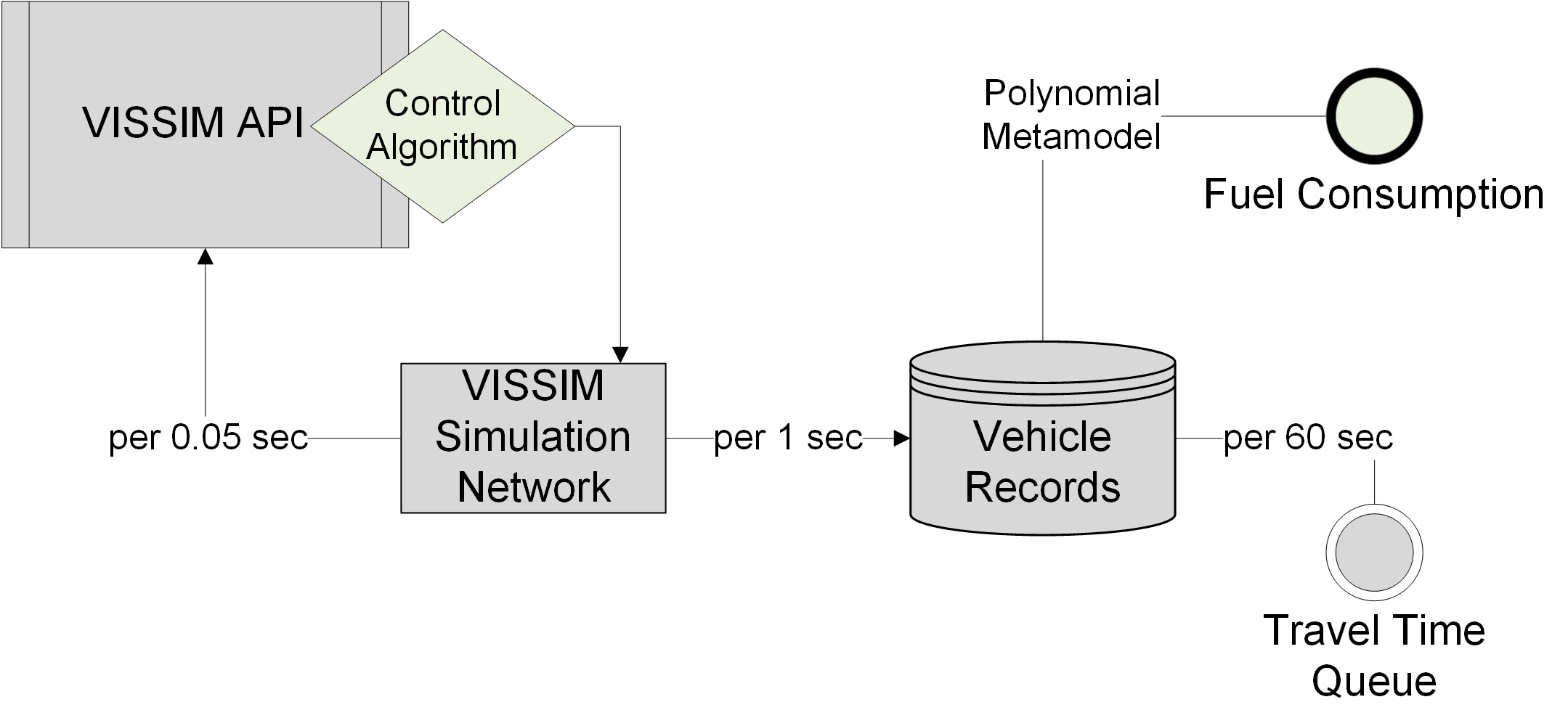}    % The printed column width is 8.4 cm.
\caption{Simulation framework in the VISSIM environment.} 
\label{fig:framework}
\end{center}
\end{figure}

To evaluate the effectiveness of the proposed approach, a simple roundabout is created and a simulation framework is established in PTV VISSIM environment (Fig.~\ref{fig:framework}). The analytical, closed-form solution described in the previous section is implemented through VISSIM Application Programming Interface (API). In this study, simulation resolution is set as 20 time steps per second. Thus, every 0.05 $sec$, with collected vehicle information (i.e., acceleration, speed, and position), the control algorithm determines the optimal acceleration/deceleration for each vehicle within the control zone, and sends back the recommended values to each vehicle under control. In the VISSIM network, vehicle trajectories are archived every 1 $sec$, and the aggregated data including travel time and delay are recorded every 60 $sec$ for network performance evaluation. Fuel consumption is estimated by using the polynomial metamodel proposed in \cite{Kamal2011} that relates vehicle fuel consumption as a function of speed $v(t)$ and acceleration $u(t)$.

\begin{figure} [!ht]
\begin{center}
\includegraphics[width=8.4cm]{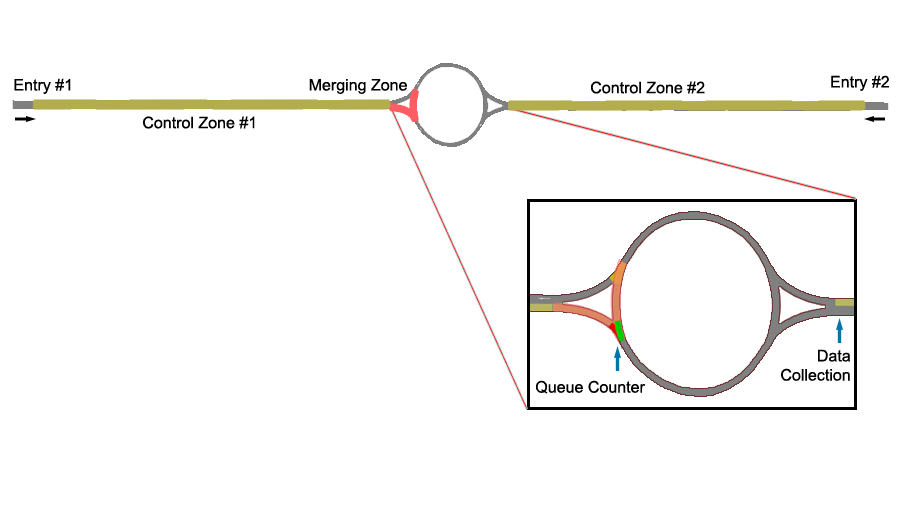}    % The printed column width is 8.4 cm.
\caption{Simulated roundabout in VISSIM.} 
\label{fig:roundaboutvissim}
\end{center}
\end{figure}

The simulation network is shown in Fig.~\ref{fig:roundaboutvissim}. The perimeter of the roundabout is approximately 200 $m$. There are two entry points from opposite approaches where vehicles will be entering the network and interacting with each other in the merging zone. In this study, two routes are defined, where eastbound traffic travel through the roundabout, and westbound traffic make a U-turn at the roundabout. Therefore, only one merging point is considered (as highlighted in light red in Fig.~\ref{fig:roundaboutvissim}), with westbound vehicles entering the roundabout and forming a circulating flow that eastbound vehicles try to merge into. The distance from each vehicle input point to the entry point of the roundabout is 320 $m$ each, including 20 $m$ entry zone. The arc length of the merging zone is set as 12 $m$, and the length of the control zone for each approach is 300 $m$ (light yellow segments in Fig.~\ref{fig:roundaboutvissim}).

According to the guidelines published by FHWA; see \cite{Robinson2000}, the recommended maximum entry speed for a roundabout varies from 6.9 $m/s$ to 13.9 $m/s$ for different roundabout designs. In this study with a single-lane roundabout, the entry speed is set as 15.6 $m/s$ for both traveling approaches and 8.9 $m/s$ is adopted as the speed limit inside the roundabout. For general-purpose vehicles, 3.0 $m/s^2$ and -4.5 $m/s^2$ are the recommended maximum acceleration and minimum deceleration settings. However, considering the presence of automated vehicle technologies, the maximum acceleration threshold here is relaxed to 4.5 $m/s^2$; see \cite{Dowling2004, Lee2013, Zhou2012}.

For 0\% penetration of CAVs (i.e., a network of non-CAVs), we apply the Wiedemann car following model; see \cite{Wiedemann1974}. In VISSIM, the minimum safe distance is defined as the distance that a driver would maintain with its leading vehicle, which is expressed as a function of standstill distance and headway time. In this study, the default 1.5 $m$ standstill distance in VISSIM is adopted and 1.2 $sec$ headway time is set for non-CAVs. 

To investigate the influence of controlled CAVs on the traffic flow, different CAV market penetration rates (MPR) are considered in the study, including 0\%, 20\%, 50\%, 80\%, and 100\% MPR. For the control of CAVs in a mixed environment, if the physically leading vehicle of a CAV is a non-CAV, the CAV will check the safety constraint continuously to make adjustment of its travel behavior. A simple on-off switch is applied in the study: the control algorithm for a CAV would be always switched on until the safety constraint is activated in terms of the distance between itself and its leading non-CAV. 

From a preliminary test, it was found that a demand of 800 vehicles per hour per lane ($vphpl$) represented a near-capacity scenario for the study network. Therefore, in this study, the input demand is set as 800 $vphpl$. For all cases, a total of 400 vehicles are dispatched from two entry points within 900 $sec$. The simulation time is set as 1200 $sec$ in order to process all the vehicles. As a result, the total travel time and fuel consumption would be the sum over all dispatched vehicles.

Five simulation runs are conducted to account for the effect of stochastic components of traffic and driver behaviors. Measures of effectiveness (MOEs), including vehicle travel time and delays are collected through VISSIM. Vehicle travel time is calculated from the entry point to the data collection point (i.e., the exit of the roundabout, as indicated in Fig.~\ref{fig:roundaboutvissim}) for each approach.

\section{Results}

\begin{figure}
\begin{center}
\includegraphics[width=8.4cm]{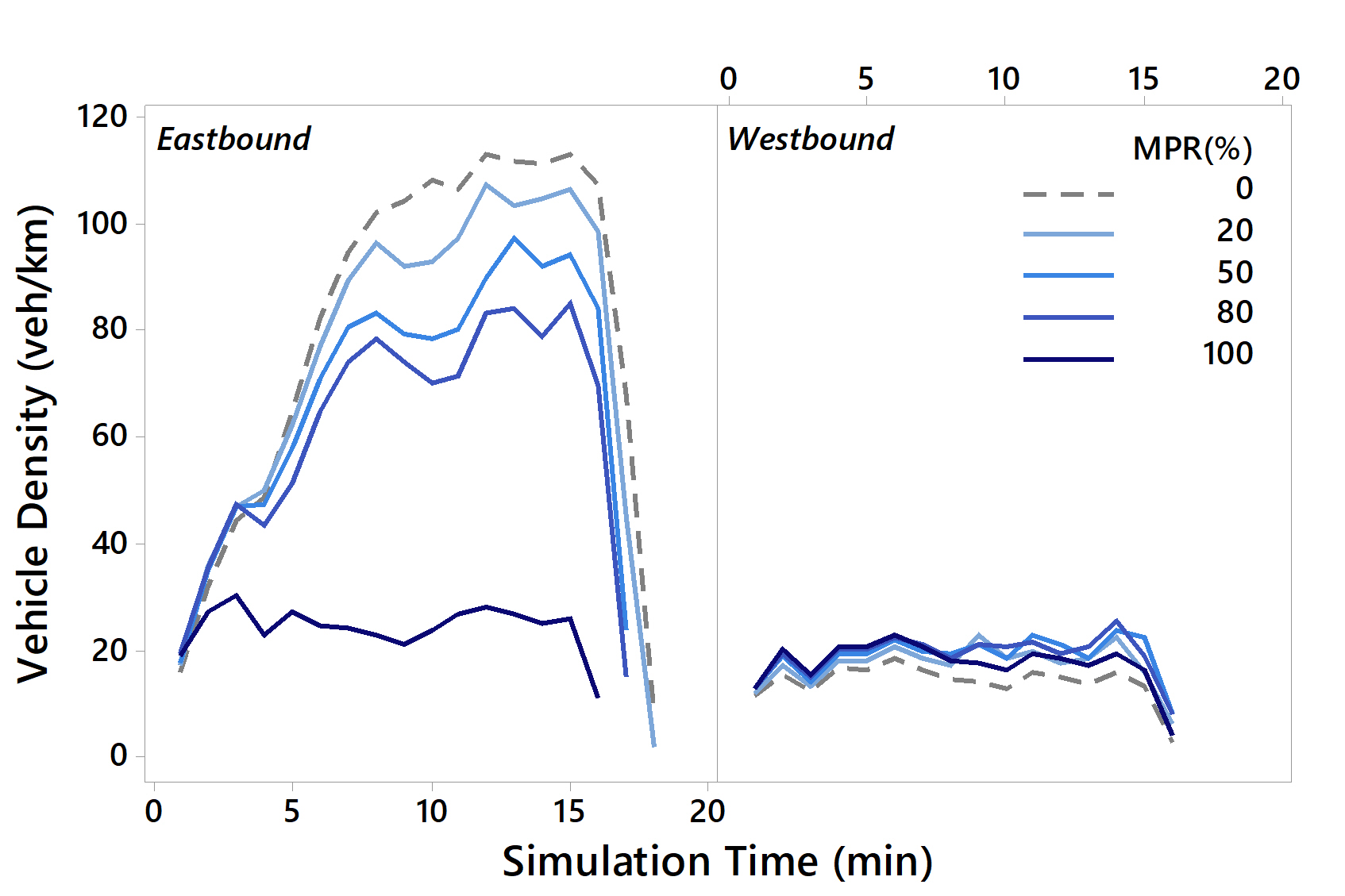}    % The printed column width is 8.4 cm.
\caption{Vehicle Density.} 
\label{fig:density}
\end{center}
\end{figure}

For the scenario of 0\% CAV penetration, the vehicles traveling eastbound have to yield to westbound traffic. If traffic density is low, the gaps between westbound traffic are generally large enough so that few eastbound vehicles need to stop in order to merge into the roundabout. However, when the demand is near capacity, it is extremely difficult for eastbound traffic to find proper gaps to merge. As a result, a queue is built up until the end of simulation. As shown in Fig.~\ref{fig:density}, for a network of non-CAVs, while the density of westbound traffic is stable during the entire simulation period, the roadway with eastbound traffic experiences increasing density as vehicles are dispatched into the network. It is not until the end of 900 $sec$ when the vehicle input stops, that the density drops and congestion is gradually released.

\begin{figure}[!ht]
\begin{center}
\includegraphics[width=8.4cm]{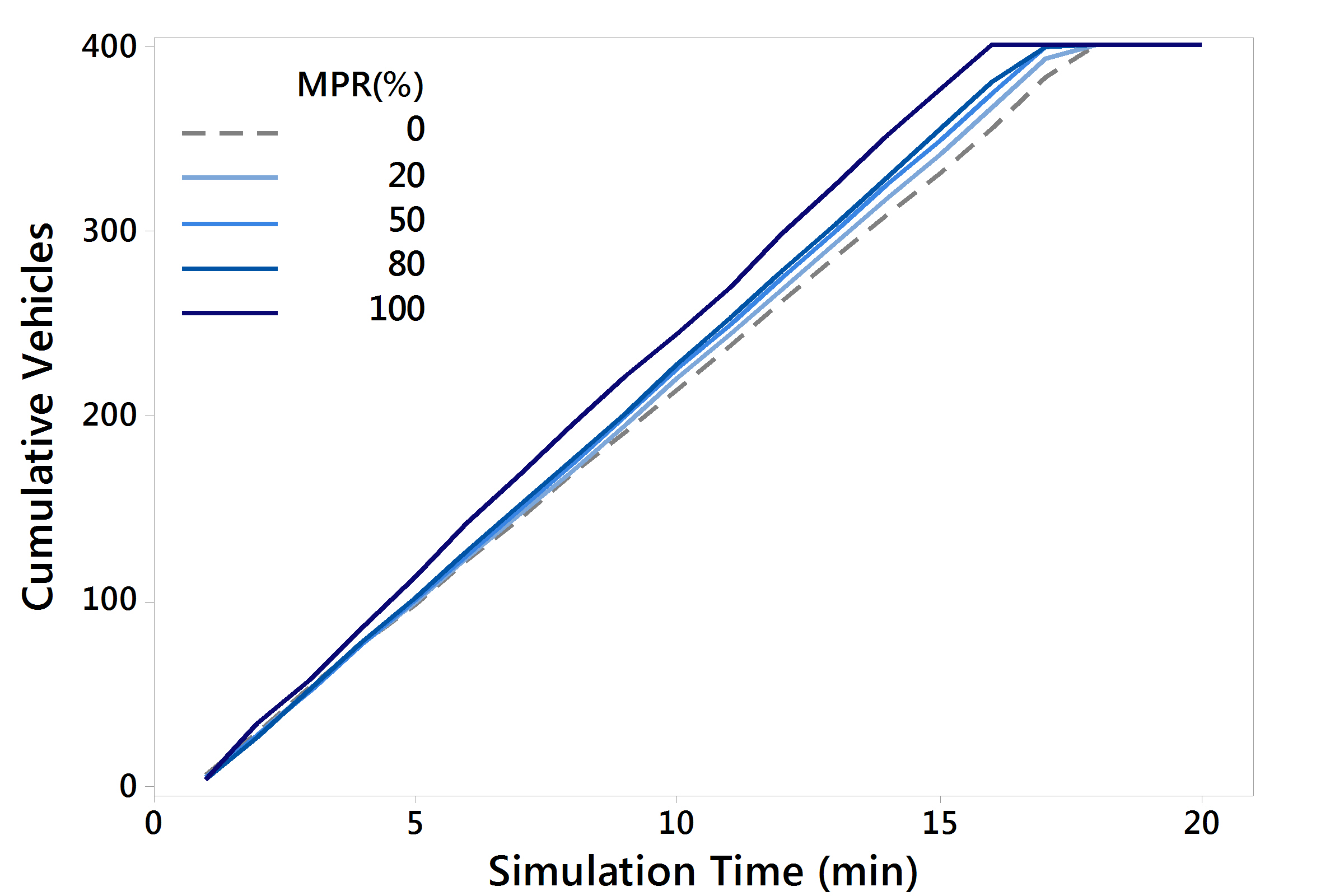}    % The printed column width is 8.4 cm.
\caption{Cumulative Vehicles.} 
\label{fig:vehs}
\end{center}
\end{figure}

It is clear that with 100\% MPR, the network performance is improved. By optimizing individual vehicle's acceleration and deceleration, eastbound vehicles are able to merge into the roundabout without stops even with high circulating flow. As shown in Fig.~\ref{fig:density}, the density of eastbound approach becomes stable (and similar to that of westbound approach) throughout the simulation period -- the network capacity is leveraged to smooth traffic flow for both traveling approaches. Therefore, during the same time period, the cumulative number of vehicles served by the roundabout increases as compared to 0\% MPR case (Fig.~\ref{fig:vehs}), leading to an improved roundabout capacity.

\begin{figure}[!ht]
\begin{center}
\includegraphics[width=8.4cm]{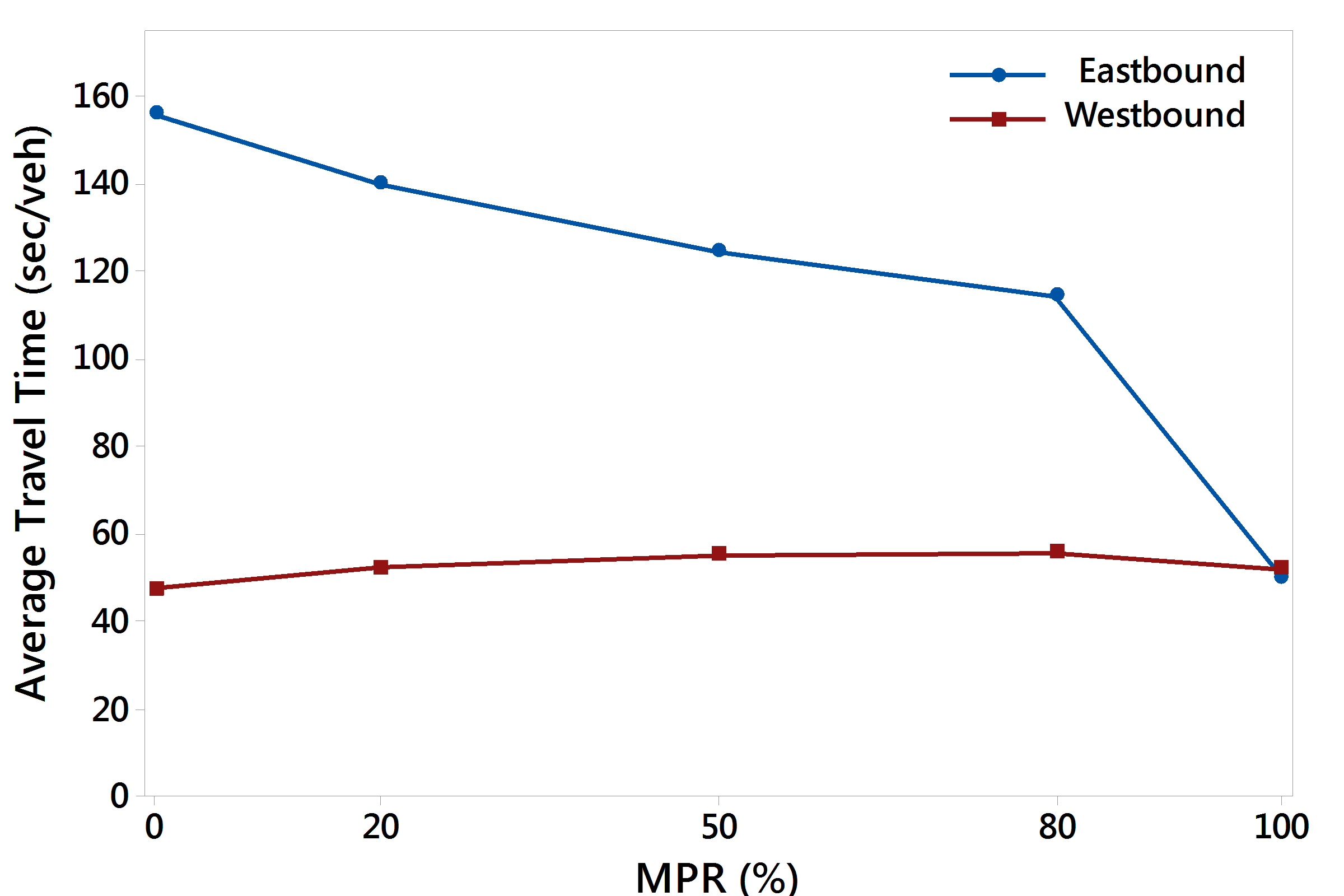}    % The printed column width is 8.4 cm.
\caption{Average travel time.} 
\label{fig:avg_travel_time}
\end{center}
\end{figure}

\begin{figure}[!ht]
\begin{center}
\includegraphics[width=8.4cm]{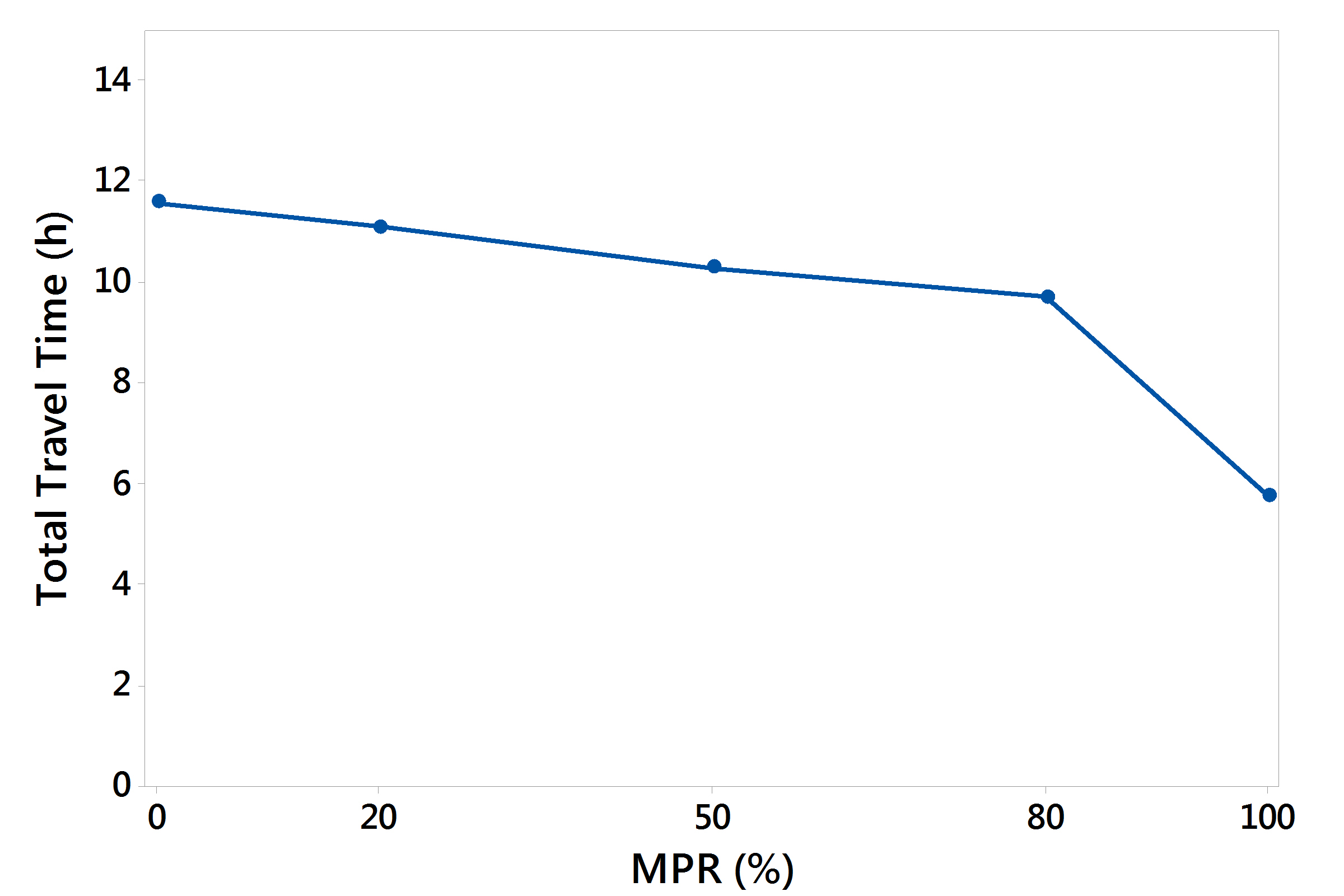}    % The printed column width is 8.4 cm.
\caption{Total network travel time.} 
\label{fig:sum_travel_time}
\end{center}
\end{figure}

\begin{figure}[!ht]
\begin{center}
\includegraphics[width=8.4cm]{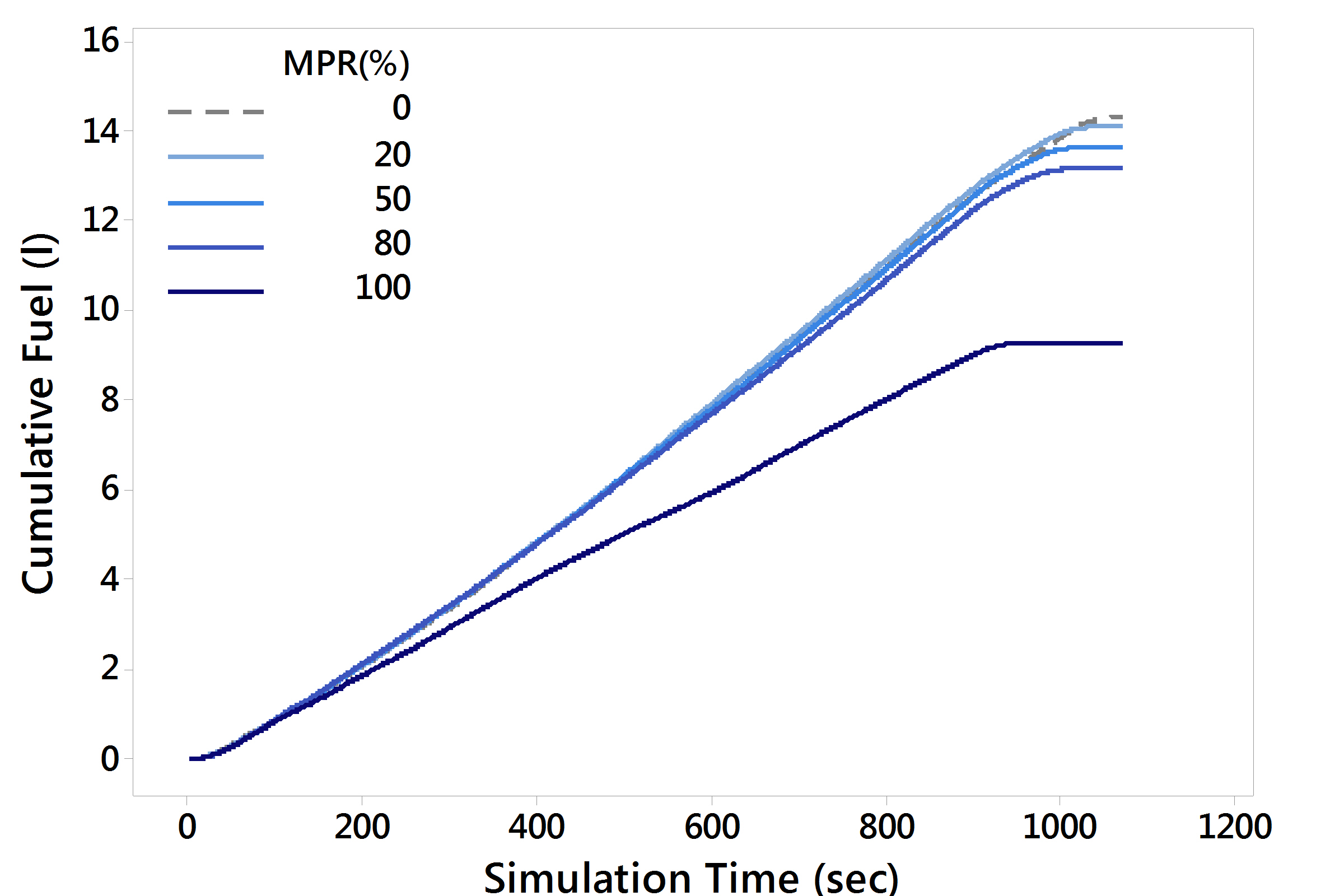}    % The printed column width is 8.4 cm.
\caption{Cumulative fuel consumption.} 
\label{fig:cum_fuel}
\end{center}
\end{figure}

In addition, through vehicle coordination under 100\% MPR, the large variation in traffic conditions for two approaches is minimized (Fig.~\ref{fig:avg_travel_time}), and the overall network travel time (Fig.~\ref{fig:sum_travel_time}) is improved significantly. As a result, a 51\% travel time savings is observed for the entire network. Furthermore, by eliminating vehicles' stop-and-go driving for eastbound traffic, transient engine operation is minimized, leading to direct fuel consumption savings as shown in Fig.~\ref{fig:cum_fuel}. 

Looking into the scenarios with mixed traffic, whereas the introduction of controlled CAVs could lead to improved network performance in terms of total network travel time and fuel consumption, the benefit is not substantial even with high MPR (Table~\ref{tb:moes}). This is discussed as follows. Under mixed traffic condition, a leading non-CAV could slow down a series of following CAVs if it stops before the roundabout to yield to the circulating flow. Especially, when the traffic demand is high, the traffic flow becomes extremely unstable and is sensitive to disturbance. In such case, even with only one stopped non-CAV, a queue can be easily formed. Furthermore, due to high circulating flow, it is hard to dissolve such queue. Therefore, under high CAV MPR (e.g., 80\%), even though there are enough CAVs traveling eastbound to form a smooth flow, the delay caused by non-CAVs that wait for safe gaps to merge could be substantial.

\begin{table}[!ht] \label{tb:moes}
\centering
\caption{Improvements in measures of effectiveness under different market penetration rates.} 
\begin{tabular}{cccc}
MPR & Travel Time & Total Delay & Fuel Consumption  \\ 
(\%) & (\%) & (\%) &(\%)\\ \hline
0 & 0 & 0 & 0 \\
20 & 4 & 15 & 1 \\
50 & 11 & 35 & 5 \\
80 & 16 & 45 & 8 \\
100 & 51 & 100 & 35 \\
\hline
\end{tabular}
\end{table}

\section{Conclusions}
The proposed approach hinges on the fact that in this new environment of massive amount of information from vehicles and infrastructure, we have intriguing and promising opportunities for enabling users to better control transportation network to reduce energy consumption, travel delays and improve safety. The efficiency of the proposed approach was investigated through a simulation environment, where a number of CAVs were controlled before they entered into a roundabout to form a smooth traffic flow. The results showed that vehicle coordination yielded significantly improved travel time and fuel consumption under 100\% CAV MPR. However, the improvement of network performance under mixed traffic condition is not as substantial as compared to 100\% MPR for a near-capacity demand scenario.

Ongoing research effort is focusing on investigating the interactions between CAVs and non-CAVs in a mixed traffic environment and improving the control algorithm to achieve better network performance with low market penetration rates of CAVs. Future research should also consider vehicle coordination for multi-lane roundabouts with increased conflicting points and weaving sections. 

\begin{ack}
This manuscript has been co-authored by UT-Battelle, LLC, under Contract No. DE-AC05-00OR22725 with the U.S. Department of Energy (DOE). The United States Government retains and the publisher, by accepting the article for publication, acknowledges that the United States Government retains a nonexclusive, paid-up, irrevocable, world-wide license to publish or reproduce the published form of this manuscript, or allow others to do so, for United States Government purposes. 
\end{ack}

\bibliography{IFAC_CTS_2018_rev3}            

\end{document}